# Gödel's Incompleteness Theorems hold vacuously

*Bhupinder Singh Anand*

*Gödel's Theorem XI essentially states that, if there is a **P**-formula [**Con**(**P**)] whose standard interpretation is equivalent to the assertion "**P** is consistent", then [**Con**(**P**)] is not **P**-provable. We argue that there is no such formula.*

**1.0 Introduction**

*Gödel's First Incompleteness Theorem*

Theorem VI of Gödel's seminal 1931 paper [Go31a], commonly referred to as "Gödel's First Incompleteness Theorem", essentially asserts:

*Meta-theorem 1*: Every omega-consistent formal system **P** of Arithmetic contains a proposition "[(**A**x)**R**(x, p)]" such that both "[(**A**x)**R**(x, p)]" and "[~(**A**x)**R**(x, p)]" are not **P**-provable.

In an earlier paper [An02b], we argue, however, that a constructive interpretation of Gödel's reasoning establishes that any formal system of Arithmetic is omega-inconsistent.

It follows from this that Gödel's Theorem VI holds vacuously.

*Gödel's Second Incompleteness Theorem*

In this paper, we now argue that Theorem XI of Gödel's paper [Go31a], commonly referred to as "Gödel's Second Incompleteness Theorem", also holds vacuously.

**1.1 Notation**

We generally follow the notation of Gödel [Go31a]. However, we use the notation "(**A**x)", whose standard interpretation is "for all *x*", to denote Gödel's special symbolism for Generalisation.

We use square brackets to indicate that the expression (including square brackets) only denotes the string[1] named within the brackets. Thus, "[(**A**x)]" is not part of the formal system **P**, and would be replaced by Gödel's special symbolism for Generalisation in order to obtain the actual string in which it occurs.

Following Gödel's definitions of well-formed formulas[2], we note that juxtaposing the string "[(**A**x)]" and the formula[3] "[**F**(x)]" is the formula "[(**A**x)**F**(x)]", juxtaposing the symbol "[~]" and the formula "[**F**]" is the formula "[~**F**]", and juxtaposing the symbol "[**v**]" between the formulas "[**F**]" and "[**G**]" is the formula "[**FvG**]".

The numerical functions and relations in the following are defined explicitly by Gödel [Go31a]. The formulas are defined implicitly by his reasoning.

**1.2 Definitions**

We take **P** to be Gödel's formal system, and define ([Go31a], Theorem VI, p24-25):

   (*i*)   "$Q(x, y)$" as Gödel's recursive numerical relation "$\sim xB(Sb(y\ 19|Z(y)))$".

   (*ii*)  "[**R**(x, y)]" as a formula that represents "$Q(x, y)$" in the formal system **P**.

---

[1] We define a "string" as any concatenation of a finite set of the primitive symbols of the formal system under consideration.

[2] We note that all "well-formed formulas" of **P** are "strings" of **P**, but all "strings" of **P** are not "well-formed formulas" of **P**.

[3] By "formula", we shall mean a "well-formed formula" as defined by Gödel.

(The existence of such a formula follows from Gödel's Theorem VII (1).)

(*iii*)   "*q*" as the Gödel-number of the formula "[**R**(x, y)]" of **P**.

(*iv*)   "*p*" as the Gödel-number of the formula "[(**A**x)][**R**(x, y)]"[4] of **P**.

(*v*)   "[p]" as the numeral that represents the natural number "*p*" in **P**.

(*vi*)   "*r*" as the Gödel-number of the formula "[**R**(x, p)]" of **P**.

(*vii*)   "17***Gen****r*" as the Gödel-number of the formula "[(**A**x)][**R**(x, p)]" of **P**.

(*viii*) "***Neg***(17***Gen****r*)" as the Gödel-number of the formula "[~][(**A**x)][**R**(x, p)]" of **P**.

(*ix*)   "$R(x, y)$" as the arithmetical interpretation[5] of the formula "[**R**(x, y)]" of **P**.

   ("$R(x, y)$" is defined by Gödel's Theorem VII ([Go31a], p29), where it is proved equivalent to "$Q(x, y)$".)

(*x*)   "***Wid***(**P**)" as the number-theoretic assertion "(E*x*)(***Form***(*x*) & ~***Bew***(*x*))".

---

[4] We note that "[(**A**x)][**R**(x, y)]" and "[(**A**x)**R**(x, y)]" denote the same formula of **P**.

[5] We define "interpretation" as in Mendelson ([Me64], p49). We note that the interpreted arithmetic expression "$F(x)$" is obtained from the formula [**F**(x)] by replacing every primitive, undefined symbol of **PA** in the formula [**F**(x)] by an interpreted arithmetic symbol (i.e. a symbol that is a shorthand notation for some uniquely meaningful, intuitive, arithmetic concept).

So the **PA**-formula [(**A**x)**F**(x)] interprets as the arithmetic expression "(A*x*)$F(x)$", and the **PA**-formula [~(**A**x)**F**(x)] as the arithmetic expression "~(A*x*)$F(x)$".

We also note that, classically, the equivalent meta-assertions "[(**A**x)**F**(x)] is a true proposition under the interpretation **I** of **PA**", and "(A*x*)$F(x)$ is a true proposition of the interpretation **I** of **PA**", are merely shorthand notations for the meta-assertion "$F(x)$ is satisfied for any given value of *x* in the domain of the interpretation **I** of **PA**" ([Me64], p51).

        ("*Wid*(**P**)" is defined by Gödel ([Go31a], p36) as equivalent to the meta-assertion "**P** is consistent".)

(*xi*)   "[**Con**(**P**)]" as the formula that represents " *Wid*(**P**)" in the formal system **P**.

(*xii*)   "*w*" as the Gödel-number of the formula "[**Con**(**P**)]" of **P** ([Go31a], p37).

(*xiii*)   "*Con*(**P**)" as the arithmetic interpretation of the formula "[**Con**(**P**)]" of **P**.

## 1.3 Gödel's semantic theses in Theorem XI

We begin by noting some semantic theses that underlie Gödel's proof of, and the conclusions he draws from, his Theorem XI ([Go31a], p36).

*Thesis 1*:   If a formula [**F**] is **P**-provable, then its standard interpretation "*F*" is a true arithmetic assertion.

*Thesis 2*:   "**P** is consistent", abbreviated "*Wid*(**P**)", can be meaningfully defined as equivalent to the number-theoretic proposition "(E$x$)(***Form***($x$) & ~***Bew***($x$))" ([Go31a], p36, footnote 63).

*Thesis 3*:   "*Wid*(**P**)" is equivalent to some arithmetic assertion "*Con*(**P**)", that is the interpretation of a **P**-formula "[**Con**(**P**)]".

## 1.4 A meta-theorem

We now argue that Gödel's *Thesis 3* is false.

*Meta-theorem 2*: There is no **P**-formula that asserts, under interpretation, that **P** is consistent.

*Proof*: We assume that Gödel's *Thesis 3* is true, and there is some formula [**Con**(**P**)] of the formal system **P** such that, under the standard interpretation:

"*Con*(**P**)" is a true arithmetical relation <==> **P** is a consistent formal system.

We take this as equivalent to the assertion "[**Con**(**P**)] is a **P**-formula that asserts, under interpretation, that **P** is consistent".

(*i*) By the definition of consistency[6], [**Con**(**P**)] is **P**-provable if **P** is inconsistent - since every formula of an inconsistent **P** is a consequence of the Axioms and Rules of Inference of **P**[7].

(*ii*) Now, if [**Con**(**P**)] were **P**-provable then, by *Thesis 1*, we would conclude, under the standard interpretation, that:

"*Con*(**P**)" is a true arithmetical assertion.

We would further conclude, by *Theses 2 and 3*, the meta-assertion:

**P** is a consistent formal system.

However, this would be a false conclusion, since **P** may be inconsistent.

(*iii*) It follows that we cannot conclude from the **P**-provability of [**Con**(**P**)] that **P** is consistent. Hence we cannot have any formula [**Con**(**P**)] such that:

"*Con*(**P**)" is a true arithmetical assertion ==> "*Wid*(**P**)" is a true number-theoretic assertion.

---

[6] We take Mendelson's Corollary 1.15 ([Me64], p37), as the classical definition of "consistency".

[7] This follows from ([Me64], p37, Corollary 1.15).

(*iv*)   We conclude that Gödel's *Thesis 3* is false, and there is no **P**-formula [**Con**(**P**)] such that:

"***Con***(**P**)" is a true arithmetical relation <==> **P** is a consistent formal system.

**1.5 A meta-lemma**

From the above we may also conclude that[8]:

*Meta-lemma 1*: Although a primitive recursive relation, and the interpretation of its formal representation, are always arithmetically equivalent, they are not always formally equivalent.

**1.6 Gödel's Proof of Theorem XI**

The question arises: Does *Meta-theorem 2* contradict Gödel's Theorem XI [Go31a]?

Now Gödel's Theorem XI [Go31a] is essentially the following assertion.

*Meta-theorem 3*: The consistency of **P** is not provable in **P**.

*Proof*:   Gödel [Go31a] argues that:

(*i*)   If **P** is assumed consistent, then the following number-theoretic assertions follow from his Theorems V, VI and his definition of "*Wid*(**P**)".

$$Wid(\mathbf{P}) => \sim Bew(17Gen\, r)$$

$$Wid(\mathbf{P}) => (Ax)\sim xB(17Gen\, r)$$

---

[8] A possible cause, and significance, of such non-equivalence is discussed in ([An01], Chapter 2), and highlighted particularly in ([An01], Chapter 2, Section 2.9).

$17Genr = Sb(p\ 19|Z(p)))$

$Wid(\mathbf{P}) \Rightarrow (Ax) \sim xB(Sb(p\ 19|Z(p)))$

$Q(x, y) \Leftrightarrow \sim xB(Sb(y\ 19|Z(y)))$

$(x)Q(x, y) \Leftrightarrow (Ax) \sim xB(Sb(y\ 19|Z(y)))$

$Wid(\mathbf{P}) \Rightarrow (Ax)Q(x, p)$

(*ii*) Assuming that [(A*x*)**R**(x, p)] asserts its own provability[9], Gödel concludes from the above that the instantiation:

*wImp*(17*Genr*),

of the recursive number-theoretic relation "*xImpy*", is a true numerical assertion.

(*iii*) From this, he concludes that:

"[**Con**(**P**)] => [(A*x*)**R**(x, p)]" is **P**-provable.

(*iv*) Now, in his Theorem VI, Gödel [Go31a] argues that, if **P** is assumed consistent, then [(A*x*)**R**(x, p)] is not **P**-provable. He thus concludes that, if **P** is assumed consistent, then [**Con**(**P**)] too is not **P**-provable.

(*v*) Implicitly assuming that *Thesis 3* is true, and so:

"*Con*(**P**)" is a true arithmetical assertion <==> "*Wid*(**P**)" is a true number-theoretic assertion,

---

[9] We argue, in a companion paper, that this semantic assumption is a consequence of construing "(A*x*)*F*(*x*)" non-constructively. It reflects the thesis that implicit Platonist assumptions underlie Gödel's reasoning both in the proofs of, and the conclusions he draws from, various meta-propositions that he asserts as Theorems in [Go31a].

Gödel further concludes that (*iv*) is equivalent to asserting that the consistency of **P** is not provable in **P**.

**1.7 Conclusion**

However, since, by *Meta-theorem 2*, Gödel's *Thesis 3* is an invalid implicit assumption, we conclude that Gödel's Theorem XI is essentially the vacuous meta-assertion:

*Meta-theorem 4*: If there is a **P**-formula [**Con**(**P**)] whose standard interpretation is equivalent to the assertion "**P** is consistent", then [**Con**(**P**)] is not **P**-provable.

(*Acknowledgement: My thanks to Professor Karlis Podnieks for suggesting the need for a comparison between Gödel's approach and that of the arguments underlying this paper.*)

(*Updated: Friday 9th May 2003 9:32:24 AM by [re@alixcomsi.com](mailto:re@alixcomsi.com)*)